\documentclass{amsart}
\usepackage{amssymb,latexsym}

\usepackage{graphicx}        
                    
\usepackage[all]{xy}

\makeindex

\newtheorem{question*}{Question}
\newtheorem{definition}{Definition}
\newtheorem{theorem}{Theorem}
\newtheorem{example}{Example}
\newtheorem{remark}{Remark}
\newcommand{\ra}{\rightarrow}
\newcommand{\N}{{\mathbb N}}

\newcommand{\Q}{{\mathbb Q}}

\newcommand{\C}{{\mathbb C}}
\newcommand{\g}{\mathfrak g}
\newcommand{\A}{\mathbb A}
\newcommand{\M}{\mathbb M}

\newcommand{\End}{\operatorname{End}}

\newcommand{\Hom}{\operatorname{Hom}}
\newcommand{\Der}{\operatorname{Der}}
\newcommand{\sgn}{\operatorname{sgn}}

\begin{document}

\title{Topics in Algebraic Deformation Theory}
\author{Anthony Giaquinto}
\address{Mathematics and Statistics, Loyola University Chicago, Chicago IL 60626
\texttt{tonyg@math.luc.edu}}
%
%

\hfill {\it{Dedicated to Murray Gerstenhaber and Jim Stasheff}}

\begin{abstract}
\noindent
We give a selective survey of topics in algebraic deformation theory
ranging from its inception to current times. Throughout, the numerous
contributions of Murray Gerstenhaber are emphasized, especially the
common themes of cohomology, infinitesimal methods, and explicit
global deformation formulas.

\noindent

\end{abstract}

\maketitle

\section{Introduction}\label{sec:1}

The concept of deformation is pervasive in mathematics. Its aim is to
study objects of some type by organizing them into continuous families
and determine how objects in the same family are related. This is the classic
problem of determining the ``moduli'' of an algebra or of a more
general structure. The moduli are, roughly, the parameters
on which the structure depends.  The idea goes back at least to
Riemann who, in his famous treatise \cite{Ri57} on abelian functions,
showed that the Riemann surfaces of genus $g$ form a single continuous
family of dimension $3g-3$. This family is the prototype of a moduli
space, a concept central to deformation theory.

 The modern era of deformations began with the pioneering work of
 Fr\"ohlicher-Nijenhuis \cite{FN57} and Kodaira-Spencer \cite{KS58} on
 deformations of complex manifolds. In particular, we see in
 \cite{FN57} the first formal use of infinitesimal (cohomological)
 methods in a deformation problem as the authors prove that if $X$ is
 a complex manifold, $T$ its sheaf of holomorphic tangent vectors, then
 there can be no perturbation of the complex structure whenever $H^1(X,T)=0$.
 In the monumental treatise
 \cite{KS58}, Kodaira and Spencer then developed a systematic theory
 of deformations of complex manifolds, including the infinitesimal and
 obstruction theories. For the case of Riemann surfaces, there are no
 obstructions as $H^2(X,T)=0$.

Algebraic deformation theory began with Gerstenhaber's seminal paper
\cite{Ge64}.  Although the analytic theory served as a model, numerous new
concepts lie within the realm of algebraic deformation theory. In
fact, all formal aspects of analytic deformations of manifolds are
special cases of those in the algebraic theory -- this will be made precise in
section \ref{diagram}.

Infinitesimal methods for algebra deformations are governed by
Hochschild cohomology. The study of infinitesimals led to the
discovery of the Gerstenhaber algebra structure on $HH^*(A,A)$, see
\cite{Ge63}. The ingredients of a Gerstenhaber algebra -- compatible
graded Lie and commutative products -- occupy a central position in
``higher structures in mathematics and physics.'' Another key higher
structure is that of the various infinity algebras: $A_{\infty},
L_{\infty}$ and their generalizations. These structures have roots in
Stasheff's landmark treatise \cite{Sta63}, which, coincidentally
appeared in the same year as \cite{Ge63}. While disjoint at the time,
the ideas in Gerstenhaber's and Stasheff's 1963 papers would become
closely intertwined in the years to come.  Indeed, Gerstenhaber called
the entire Hochschild cohomology $HH^*(A,A)$ the ``infinitesimal
ring'' of $A$ in \cite{Ge63}, even though only the components of
$HH^i(A,A)$ with $i\leq 3$ had natural interpretations related to
infinitesimals and obstructions. But more than 30 years later, it became
well known that the entire Hochschild cohomology $HH^*(A,A)$
\textit{is} the space of infinitesimals of deformations of $A$ to
an $A_{\infty}$ algebra.

This survey only represents a sampling of ideas in algebraic
deformation theory. None of the discussion is new, except for the
results in the last section on algebra variations. Many ideas are only
sketched and proofs are omitted. More topics are left out than
included. In particular, the theory of deformation quantization is
largely left out -- the reader is referred to Sternheimer's
contribution \cite{Ste} in this volume and the references therein for
the important physical perspective. Also left out is the deformation
theory of infinity algebras, gerbes, stacks, chiral algebras, affine
and dynamical quantum groups, and the like. More comprehensive surveys
of algebraic deformation theory and quantization can be found in the
excellent treatises \cite{GS88} \cite{CKTB05} and \cite{DMZ07}.

\section*{Algebraic Deformations}

Let $A$ be an associative algebra over a commutative ring $k$.

\begin{definition} A \emph{formal deformation} of $A$ is a $k[[t]]$-algebra $A_t$ which is flat
and $t$-adically complete as a $k[[t]]$-module, together with an isomorphism $A\simeq A_t/tA_t$.
\end{definition}
For every deformation, there is a $k[[t]]$-module isomorphism between
$A_t$ and $A[[t]]$. Once such an isomorphism is fixed, the
multiplication in $A_t$ is necessarily of the form
$$\mu_t: A[[t]]\otimes_{k[[t]]}A[[t]]\ra A[[t]] \quad \mbox{with}\quad
  \mu_t(a,b)=ab + \mu_1(a,b)t +\mu_2(a,b)t^2 + \cdots$$ where $ab$
  represents the multiplication of $A$ and each $\mu_i\in \Hom
  (A\otimes A,A)$ is extended to be $k[[t]]$-bilinear. Setting
  $\mu_0(a,b)=ab$, we have $\mu_t=\sum \mu_it^i$ and $\mu_t(a,b)$ will be
  denoted $=a*b$.

It is clear that one can consider formal deformations of other
algebraic structures (Lie algebras, bialgebras, algebra homomorphisms,
etc.\ ) by modifying the above definition to suit the appropriate
category. The deformation of an algebraic structure is usually
subjected to the same equational constraints as the original structure:
\begin{itemize}
\item Associative algebra $A$: $(a*b)*c = a*(b*c)$. If $A$ is
  commutative we can also require $a*b=b*a$.
\item Lie algebra $L$: $[a,b]_{t} = -[b,a]_{t}$ and the Jacobi
  identity for $[a,b]_{t}=[a,b]+[a,b]_1t+[a,b]_2t^2+\cdots.$
\item Bialgebra $B$: associativity of $*$, coassociativity of
  $\Delta_t(a)=\Delta(a) + \Delta_1(a)t +\cdots$, and
  $\Delta_t(a*b)=\Delta_t(a)*\Delta_t(b).$
\item Algebra homomorphism $\phi:A\rightarrow A'$: $\phi_t(ab)
  =\phi_t(a)\phi_t(b)$ with
  $\phi_t(a)=\phi(a)+\phi_1(a)t+\phi_2(a)t^2+\cdots .$
\end{itemize}

Even though in this note we are concerned with formal deformations,
there are many important and explicit instances for which the deformed
products converge or are even polynomial in $t$ when $k=\mathbb R$ or
$\mathbb C$.

\section{The deformation philosophy of Gerstenhaber}\label{GePhil}

The pioneering principle of Gerstenhaber is that the equational
constraints above can be naturally interpreted in terms of the
appropriate cohomology groups and higher structures on them.  In
particular, the infinitesimal (linear term of the deformation) is a
cocycle in the cohomology group -- Hochschild $HH^2(A,A)$ in the
associative case, Harrison $Har^2(A,A)$ in the commutative case,
Chevalley-Eilenberg $H^2_{CE}(L,L)$ in the Lie case,
Gerstenhaber-Schack $H^2_{GS}(B,B)$ in the bialgebra case, and the
diagram cohomology $H_d^2(\phi,\phi)$ in the algebra homomorphism
case. Moreover, the obstructions to extending infinitesimal and $n$-th
order deformations to global ones are controlled by the differential
graded Lie algebra structure on the cohomology.

In the associative case, the graded Lie structure (and much more) was
laid out in \cite{Ge63}. There it was shown that the Hochschild
cohomology $HH^*(A,A)=\bigoplus_{n\geq0}H^n(A,A)$ has a remarkably rich
structure consisting of two products,
\begin{itemize}
\item A graded commutative product where $\deg HH^p(A,A) = p$,
\item A graded Lie product where $\deg HH^p(A,A)=p-1$,
\item $[\alpha, -]$ is a graded derivation of the commutative product.
\end{itemize}
A graded $k$-module satisfying the above conditions is a {\em
  {Gerstenhaber algebra}}. Other notable examples are $\bigwedge^* L$
(where $L$ is a Lie algebra), $H^*(X,\bigwedge ^*T)$ (where $X$ is a
manifold and $T$ is its sheaf of tangent vectors), and the diagram
cohomology $H_d^*(\mathbb A, \mathbb A)$ of an arbitrary presheaf $\A$
of $k$-algebras (to be defined in Section \ref{diagram}). The
Chevalley-Eilenberg, Harrison, and bialgebra cohomology cohomologies
carry graded Lie brackets, but are not Gerstenhaber algebras in
general.

In \cite{Ge63}, the commutative and Lie products on the Hochschild
cohomology $HH^*(A,A)$ are defined at the cochain level and are proved
to descend to the level of cohomology. An intrinsic interpretation of
the graded Lie structure was given by Stasheff in \cite{Sta93}. There
he proved that the Gerstenhaber bracket coincides with the natural
graded bracket on $Coder(BA, BA)$, where $BA$ is the bar complex of
$A$.

Returning to the equational constraints for a deformation of an
algebra $A$, the associativity of $\mu_t$ can succinctly be expressed
in terms of the Gerstenhaber bracket as $[\mu_t,\mu_t]=0$. Writing
$\mu_t= \mu_0+\mu'$ it follows that
$2[\mu_0,\mu']+[\mu',\mu']=0$. Since the coboundary in the shifted
Hochschild complex $C^*(A,A)[1]$ is $\delta =[\mu_0,-]$, the first
summand is $2\delta \mu'$. Thus we arrive at the fundamental
associativity equivalences
\begin{equation}\label{MC}
  \mu_t\quad\text{associative}\quad \Longleftrightarrow \quad [\mu_t,\mu_t]=0\quad  \Longleftrightarrow \quad \delta(\mu')+\frac{1}{2}[\mu',\mu']=0.
\end{equation}
Thus $\mu_t$ is associative if and only if $\mu'$ satisfies the
Maurer-Cartan equation.  Although not explicitly stated as such, the
idea that deformations are governed by a differential graded Lie
algebra and solutions to the Maurer-Cartan equation goes back to
Gerstenhaber's original paper \cite{Ge64}.

\section{Algebras with Deformations}

The search for deformations of an algebraic structure $A$ begins with
the appropriate cohomology group (usually $H^2(A,A)$) which comprises
the infinitesimals. Given an infinitesimal $\mu_1$, the basic question
is whether it can be integrated to a full deformation or not. In other
words, is it possible to find $\mu_2, \mu_3,\ldots$ such that
$\mu'=\sum_{i\geq1} \mu_it^i$ satisfies the Maurer-Cartan equation. Of
course, the vanishing of the obstruction group (usually $H^3(A,A)$)
guarantees that any infinitesimal is integrable, but this is rarely
the case. A necessary condition in general is that the \emph{primary
  obstruction}, $[\mu_1,\mu_1]$, must equal zero. There are then
higher obstructions which must vanish in order for $\mu_1$ to be
integrable.  Thus, a reasonable starting point for deformations is to
first determine which infinitesimals have a vanishing primary
obstruction. Remarkably, in several fundamentally important cases, all
infinitesimals $\mu_1$ with $[\mu_1,\mu_1]=0$ can be integrated. In
fact, some of the most celebrated theorems in deformation theory are
expressions of this phenomenon.

Perhaps the most studied algebra deformations are those which lie in
the realm of deformation quantization, a concept introduced in the
seminal paper \cite{BFFLS78}. Suppose $X$ is a real manifold and
$A=C^{\infty}(X)$. Then $HH^2(A,A)$ can be identified with the space
of bivector fields $\alpha \in \Gamma (X,\wedge^2 T)$, and the primary
obstruction to $\alpha$ is the Schouten bracket $[\alpha,
  \alpha]$. The condition $[\alpha, \alpha]=0$ asserts that $\alpha$
determines Poisson structure on $X$. In \cite{BFFLS78} it was asked
whether any Poisson structure can be quantized. The affirmative answer
to this question is one of the jewels of deformation theory.

\begin{theorem}[Kontsevich]\label{Kontsevich} Any Poisson manifold can be quantized.
  More generally, there is, up to equivalence, a canonical
  correspondence between associative deformations of the algebra $A$ and formal
  Poisson structures $\alpha_t=\alpha_1t+\alpha_2t^2+\cdots$ on $A$.
\end{theorem}

Also in \cite{Ko97} is a remarkable explicit quantization formula for
$X=\mathbb R^n$. The formula involves certain weighted graphs which
determine the $*$-product expansion. A physical interpretation of the
deformation quantization formula in terms of path integrals of models
in string theory was made precise by Cattaneo and Felder in
\cite{CF00}. In the case where $X$ is a smooth algebraic variety,
quantization of Poisson brackets is also possible, but significant
modifications of the approach are necessary, see \cite{Ko01},
\cite{VdB07}, and \cite{Ye05}.

Another case where the primary obstruction to integrating an
infinitesimal is the only one is in the realm of quantum
groups. Consider a Lie bialgebra $\mathfrak a$. The cocommutator,
$\delta:\mathfrak a\rightarrow \mathfrak a \otimes \mathfrak a$, can
be extended to an infinitesimal deformation of the coalgebra structure
of $U\mathfrak a$ which is compatible with the Lie bracket on
$\mathfrak a$, and hence the multiplication in $U\mathfrak a$. The
cocommutator may thus be viewed as an infinitesimal whose primary
obstruction vanishes. Drinfel'd asked in \cite{Dr92} whether any Lie
bialgebra can be quantized.  The affirmative answer to this question
is another famous result in deformation theory.

\begin{theorem}[Etingof-Kazhdan]\label{Etingof-Kazhdan} Any Lie bialgebra can be quantized.
  That is, if $\mathfrak a$ is a Lie bialgebra, then there exists a
  Hopf algebra deformation of $U\mathfrak a$ whose infinitesimal is
  the cocommutator of $\mathfrak a$.
\end{theorem}
The quantization of $U\mathfrak a$ depends on a choice of Drinfel'd
associator. It is known that associators are not unique and are
notoriously difficult to compute with.

Many of the results pertaining to quantization of solutions to the
various types of classical Yang-Baxter equation can also be viewed as
examples of the phenomenon that in certain situations, the primary
obstruction to integrating an infinitesimal structure is the only
one. Some of these instances will be discussed in Section
\ref{bialgebra}.

In general, the condition $[\mu_1,\mu_1]=0$ does not guarantee
that an infinitesimal $\mu_1$ is integrable. The earliest known
example is geometric in nature and predates the algebraic theory. In
\cite{Do60}, Douady exhibited an example of an infinitesimal
deformation (in the Kodaira-Spencer sense) of the Heisenberg group
whose primary obstruction vanishes, yet its secondary obstruction, a Lie-Massey
bracket fails to vanish. More
recently and in the algebraic case, Mathieu has given examples of
commutative Poisson algebras which cannot be quantized, see
\cite{Mat97}.

\section{Algebras without deformations}\label{rigid}
A deformation $A_t$ of an algebra $A$ is {{\it trivial}} if there is a
  $k[[t]]$-algebra isomorphism $\Phi_t:A_t\rightarrow A[[t]]$ which reduces
  the identity modulo $t$.
An algebra is
\emph{rigid} if it has
no non-trivial deformations. The cohomology results of Section
\ref{GePhil} provide the first elementary result in deformation
theory.

\begin{theorem}\label{rigid-theorem} If $H^2(A,A)=0$ then $A$ is rigid.
\end{theorem}

Algebras which satisfy the hypothesis of Theorem \ref{rigid-theorem} are called
\textit{absolutely rigid}. Here are some notable examples of absolutely rigid
algebras in various categories.
\begin{itemize}

\item Any separable $k$-algebra  $A$ is rigid as these are characterized by  $HH^n(A,-)$ for all $n\geq 1$.

%

%

\item The enveloping algebra $U\mathfrak g$ of a finite
  dimensional semisimple Lie algebra $\mathfrak g$ is rigid as an
  algebra as $HH^n(U\mathfrak g,U\mathfrak g)=0$ for $n\geq 1$. It
  does admit deformations as a Hopf algebra.
\item The coordinate ring $\mathcal O(V)$ of a smooth affine variety
  $V$ is rigid as a commutative algebra as $Har^n(\mathcal
  O(V),\mathcal O(V))=0$ for $n\geq 1$ -- a more precise interpretation of
  the cohomology $\mathcal O(V)$ will be given in Section
  \ref{Hodge}. It does admit noncommutative deformations however.
\item The $m$-th Weyl (Heisenberg) algebra $A_m$ is rigid as $HH^n(A_m,A_m)=0$ for $n\geq 1$.

\end{itemize}

The converse of Theorem \ref{rigid-theorem} is known to be false in many instances.
Richardson has provided in \cite{Ri67}
examples of rigid Lie algebras $L$ with $H^2_{CE}(L,L)\neq 0$. In the
associative case, Gerstenhaber and Schack have given examples of rigid
associative algebras when $\operatorname{char}(k)=p$ in \cite{GS86}. Remarkably,
the question of whether there exist rigid associative algebras $A$ in
characteristic zero with $HH^2(A,A)\neq 0$ is still an open question
even in the case where $A$ is a finite-dimensional $\mathbb C$-algebra.

The rigid algebras of Gerstenhaber and Schack in
$\operatorname{char}(k)=p$ are not everyday examples. The smallest
rigid algebra constructed with $HH^2(A,A)\neq 0$ is a
$669$-dimensional algebra over $\mathbb F_2$. The algebra is a poset
algebra of a suspension of a triangulation of the projective
plane.

The proof of the rigidity of these algebras despite non-zero
$HH^2(A,A)$ is based on an elementary but fundamental theorem
of \cite{GS86} concerning
\textit{relative} Hochschild cohomology. If $S$ is a subalgebra of $A$
then a cochain $F\in C^n(A,A)$ is \textit{S-relative} if
\begin{multline}F(sa_1,\ldots,a_n)=sF(a_1,\ldots,a_n),\quad
  F(a_1,\ldots,a_ns)=F(a_1,\ldots,a_n)s,\\ \text{and}\quad F(\ldots,
  a_is,a_{i+1},\ldots)=F(\ldots, a_i,sa_{i+1},\ldots)\end{multline}
for all $s\in S$ and $a_i\in A$. Further, an $S$-relative cochain $F$ is
\textit{reduced} if $F(a_1,\ldots,a_n)=0$ whenever any $a_i\in S$.

\begin{theorem}\label{relative} If $S$ is a
  separable subalgebra of $k$, then the inclusion of the complex of
  reduced $S$-relative cochains into the full Hochschild complex
  induces an isomorphism of cohomology.
\end{theorem}

The theorem significantly reduces cohomology computations whenever $A$
has a large separable subalgebra. For example, a poset algebra is a
subset of the algebra of $n\times n$ matrices and one may take $S$ to
be the diagonal matrices. Using Theorem \ref{relative} it is
elementary to show that the Hochschild cohomology of the poset algebra
coincides with the simplicial cohomology of the geometric realization
or nerve (see Section \ref{diagram}) of the poset.

Another application of the theorem is in the computation of the
cohomology of a crossed product algebras. Let $V$ be a finite
dimensional $k$-vector space, $SV$ its symmetric algebra, and $G$ a
finite group which acts on $V$, and hence on $SV$. For $x\in SV$ and
$g\in G$ denote the action of $g$ on $x$ by $x^g$. The crossed product
$SV\ltimes G$ has underlying space $SV\otimes kG$, with the usual
multiplication in $SV$ and $kG$, and relations $(1\otimes g)(x\otimes
1) = x^g\otimes g$ for $x\in SV$ and $g\in G$. For simplicity, we omit
the tensor product symbol and write an element $x\otimes g$ simply as
$xg$.

When $|G|$ is invertible in $k$, then Maschke's Theorem asserts that
$kG$ is separable -- let us assume this here. In this case, we may
compute $HH^*(SV\ltimes G,SV\ltimes G)$ using reduced $kG$-relative
cochains. For such a cochain $F$, we have
$$F(x_1g_1,x_2g_2 \ldots, x_ng_n)=F(x_1,x_2^{g_1}\, \ldots \,
,x_n^{g_1\cdots g_{n-1}})\,g_1\cdots g_n.$$ The right side of the
equation is an element of $C^n(SV, SV\ltimes kG)$ and it is easy to
see that it is $G$-invariant in the sense that
$$gF(x_1,\ldots ,x_n)g^{-1} = F(gx_1g^{-1},\ldots, g x_n g^{-1}).$$
Therefore we obtain
$$HH^*(SV\ltimes G,SV\ltimes G)\simeq (HH^*(SV,SV\ltimes G))^G.$$
A complete computation of the cohomology can be found in many sources, and one
which explicitly uses $kG$-relative cohomology is \cite{Pi06}.
There is interest in crossed product cohomology and deformations as
they have geometric implications for orbifolds (see \cite{CGW04}) and
in the theory of symplectic reflection algebras (see \cite{EG02} and
Section \ref{relations} of this survey).

It is interesting to note that some of the rigid algebras listed above
naturally appear in parametrized families, a seeming contradiction to
the general theory. For example, consider $\mathcal
O_{\lambda}=\C[x,y]$ with $y^2=x(x-1)(x-\lambda)$, the ring of regular
functions on an affine elliptic curve. If $\lambda$ is close to
$\lambda'$ then $\mathcal O_{\lambda}\ncong \mathcal
O_{\lambda'}$. Nevertheless, $HH^2(\mathcal O_{\lambda},\mathcal
O_{\lambda})=0$. According to Kontsevich, the problem here is that the
variety, being affine, is not compact and formal deformation theory
for non-compact objects can give ``nonsensical'' results, see
\cite{Ko01}. In the associative case, a ``compact'' object is a
finite-dimensional algebra and so we expect other nonsensical results
for some infinite-dimensional algebras. Here is such an example:
the first quantized Weyl algebra
$A_q=\C \langle x,y\rangle / (qxy-yx-1)$ is not isomorphic to $A_1$ for $q$ near
1. However, as noted above, $HH^2(A_1,A_1)=0$. To put this into a
formal deformation theoretic perspective, let $q=1+t$. Then there is
indeed an analytic isomorphism $\phi: A_1[[t]]\rightarrow A_q[[t]]$,
but it has zero radius of convergence. A similar phenomenon happens
for the situation with $O_{\lambda}$ and $O_{\lambda'}$. The problem
is that passing to the formal power series versions of these algebras
has trivialized the deformations.

The above examples suggest that the classic deformation theory of a
single algebra does not always detect the dependency of an algebra on
parameters. However, the more general \textit{diagram} cohomology
theory of Section \ref{diagram} can detect such dependencies, but does
not show how the algebras vary with the parameters. The construction
of the algebras with varying moduli can be accomplished through the idea of a
\textit{variation} of algebras. This concept will be addressed in Section
\ref{variation}.

\section{Universal Deformation Formulas}\label{UDF}

The process of constructing deformations using the infinite
step-by-step procedure of extending deformations of order $n$ to $n+1$
for each $n\geq 1$ is impractical.  There are instances though in
which a closed form for $\mu_t$ is known. One is the explicit
quantization of Poisson brackets on $\mathbb R^n$ given in
\cite{Ko97}.  Another comes from the use of ``universal deformation
formulas'' which are, in essence, Drinfel'd twists which act on a
certain classes of algebras. The prototypical example of this type of
formula was given by Gerstenhaber in \cite{Ge68}. There it was
observed that if $\phi$ and $\psi$ are commuting derivations of any
associative algebra (in characteristic zero), then $a*b=\displaystyle
\sum \phi^n(a)\psi^n(b)\frac{t^n}{n!}$ is associative. The most
famous use of this idea gives the Moyal product. For example, if
$A=k[x,y]$ with $\phi=\partial_x$ and $\psi=\partial_y$, then we have
$x*y-y*x=t$, and the deformation is isomorphic to the first Weyl
algebra as long as $t\neq 0$. When $\phi=x\partial_x$ and
$\psi=y\partial y$ then the deformation is graded and isomorphic to
the skew-polynomial ring $k\langle x,y\rangle/(qxy-yx)$ with $q=e^t$. These
examples can of course be extended to higher dimensions.

\begin{definition} Suppose $B=(B,\Delta,1,\epsilon)$ is a bialgebra with comultiplication
$\Delta$, unit $1$, and counit $\epsilon$. A universal deformation
  formula (UDF) based on $B$ is an element $F\in (B\otimes B)[[t]]$
  such that
$$((\Delta \otimes 1)(F))(F\otimes 1) = ((1\otimes
  \Delta)(F))(1\otimes F)\quad\text{and}\quad (\epsilon \otimes
  1)F=(1\otimes \epsilon)F=1\otimes 1.$$
\end{definition}
The virtue of a UDF is that for any $B$-module algebra $A$, the
product $a*b=\mu\circ F(a\otimes b)$ is associative and hence is a
deformation of $A$, see \cite{GZ98}.

\begin{example}\label{MoyalUDF}
Suppose $B$ is commutative and let $r\in P\otimes P$ where $P$ is the
space of primitive elements. Then $F=\exp(tr)$ is a UDF. Primitive
elements of $B$ act as derivations of any $B$-module algebra and so
this UDF gives a wide range of Moyal-type deformations.
\end{example}

\begin{example}\label{S2UDF} Let $B=U\mathfrak s$, where $\mathfrak s$ is the Lie algebra
with basis $\{H,E\}$ and relation $[H,E]=E$. Set $H^{\langle
  n\rangle}= H(H+1)\cdots (H+n-1)$. Then $\displaystyle
F=\sum\frac{t^n}{n!}H^{\langle n\rangle}\otimes E^n$ is UDF. For an
example of its use, take $A=k[x,y]$ with the derivations
$H=x\partial_x$ and $E=x\partial_y$. The deformed algebra has the relation
$x*y-y*x=tx^2$ and is the Jordan quantum plane. Numerous generations
of this UDF can be found in \cite{KLO01}, \cite{LS02} and the
references therein.
\end{example}

\begin{example}\label{twist} Let $\g \otimes \g$ be a Lie algebra and let
$r\in \g \wedge \g$ satisfy $[r,r]=0$, where $[-,-]$ is
  the Schouten bracket. Drinfel'd has shown in \cite{Dr83} that there
  exists a UDF $F=1\otimes 1+tr+O(t^2)$. Examples \ref{MoyalUDF} and
  \ref{S2UDF} are of this form. It should be noted that $[r,r]=0$
  means that $r$ is a skew-symmetric solution of the classical Yang-Baxter
  equation.
\end{example}

\begin{example} Let $B$ be the bialgebra generated by $\{D_1, D_2, \sigma\}$ with
relations $D_1D_2=D_2D_1, \quad D_i\sigma = q\sigma D_i \,\, (i=1,2)$, and
comultiplication

$$\Delta(D_1)=D_1\otimes \sigma+1\otimes D_1,\quad
  \Delta(D_2)=D_2\otimes 1+\sigma \otimes D_1,\quad
  \Delta(\sigma)=\sigma\otimes \sigma.$$
Then $F=\exp_q(tD_1\otimes D_2)$ is a UDF, where the $q$-exponential
 is the usual exponential series with $n!$ replaced by $n_q!$.

 Note that for any $B$-module algebra, $\sigma$ acts as an
 automorphism and. the elements $D_1, D_2$ act as commuting skew
 derivations with respect to $\sigma$. Thus, this UDF provides
 $q$-Moyal type deformations. For example, it can be used to deform
 the quantum plane $k\langle x,y\rangle/(qxy-yx)$ to the first
 quantized Weyl algebra $A_q=k\langle x,y\rangle/(qxy-yx-1)$. Formulas of
 this type were also used in \cite{CGW04} to deform certain crossed
 products $SV\ltimes G$.

\end{example}

Recently, universal deformation formulas have arisen naturally in the
work of Connes and Moscovici on Rankin-Cohen brackets and the Hopf
algebra $H_1$ of transverse geometry, see \cite{CM04}. Rankin-Cohen
brackets are families of bi-differential operators on modular
forms. These brackets can be assembled to give universal deformation
formulas. Some applications appear in \cite{CM04}. The formulas based on $H_1$
are also connected to certain topics in deformation quantization as it
relates to the Poisson geometry of groupoids and foliations. see
\cite{BTY07}.

\section{Commutative Algebras and the Hodge Decomposition}\label{Hodge}

Let $A$ be a commutative algebra over a field of characteristic
zero. In \cite{Ba68}, Barr proved that the Harrison cohomology
$Har^n(A,A)$ is a direct summand of the Hochschild cohomology
$HH^n(A,A)$. The key to this splitting was Barr's discovery of an
idempotent $e_n$ in $\Q S_n$, the rational group algebra of the
symmetric group. The symmetric group acts on $C^n(A,M)$ (the
Hochschild $n$-cochains of $A$ with coefficients in a symmetric
$A$-bimodule $M$) via $\sigma F(a_1,\dots, a_n) = F(a_{\sigma
  1},\dots, a_{\sigma n})$. Barr proved that $\delta
(e_nF)=e_{n+1}(\delta F)$, where $\delta$ is the Hochschild coboundary
operator. Thus $HH^n(A,M)$ splits as $e_nHH^n(A,M)\oplus
(1-e_n)HH^n(A,M)$, and the latter piece is $Har^n(A,M)$. Barr's work
received little attention until 1987 when Gerstenhaber and Schack
extended the splitting, see \cite{GS87}. In $\Q S_n$ there are $n$
mutually orthogonal idempotents $e_n(1),\dots,e_n(n)$ with the
property that $\delta (e_n(r)F) = e_{n+1}(r)(\delta F)$ for all $F\in
C^n(A,M)$. The relation between the idempotents and coboundary give
the following fundamental theorem.

\begin{theorem}[Hodge Decomposition] Suppose
  $\operatorname{char}(k)=0$ and let $A$ be a commutative algebra and
  $M$ a symmetric $A$-bimodule. Then there is a splitting
$$HH^n(A,M)=HH^{1,n-1}(A,M)\oplus HH^{2,n-2}(A,M)\oplus\cdots \oplus
HH^{n,0}(A,M)$$ where $HH^{r,n-r}(A,M)$ is the cohomology of the complex
$e_{*}(r)C^*(A,M)$. \end{theorem}

Around the same time as the Gerstenhaber-Schack
paper \cite{GS87}, Loday, using different techniques, exhibited a
splitting of the Hochschild and cyclic cohomologies of a commutative
algebra, see \cite{Lo88}, \cite{Lo89}.

The idempotents $e_n(r)$, which have independent interest apart from
cohomology, are most easily described using the following elegant
generating function discovered by Garsia in \cite{Ga90}:

$$\sum e_{n}^{(r)}x^r =\frac{1}{n!}\sum_{\sigma\in
  S_n}\sgn(\sigma)(x-d_{\sigma})(x-d_{\sigma}+1)\cdots
(x-d_{\sigma}+n-1)\, \sigma$$ where $d_{\sigma}$ is the number of
descents in $\sigma$, i.e. the number of $i$ with
$\sigma(i)>\sigma(i+1)$.

The following diagram, in which $HH^{i,n-i}(A,M)$ is abbreviated as
$H^{i,n-i}$, is instructive in understanding the Hodge decomposition.

$\begin{array}{ccccccccccc}
&&&&&&&&&&\\
&&&&&&&&&&\makeatletter
\def\Ddots{\mathinner{\mkern1mu\raise\p@
\vbox{\kern7\p@\hbox{.}}\mkern2mu
\raise4\p@\hbox{.}\mkern2mu\raise7\p@\hbox{.}\mkern1mu}}
\makeatother \Ddots\\
&&&&&&&&&&\\
&&&&&&&&\boxed{H^{5,0}}&\longrightarrow&\cdots\\
&&&&&&&&&&\\
&&&&&&\boxed{H^{4,0}}&\longrightarrow&\boxed{H^{4,1}}&\longrightarrow &\cdots\\
&&&&&&&&&&\\
&&&&\boxed{H^{3,0}}&\longrightarrow&\boxed{H^{3,1}}&\longrightarrow&\boxed{H^{3,2}}&\longrightarrow &\cdots\\
&&&&&&&&&&\\
&&\boxed{H^{2,0}}&\longrightarrow &\boxed{H^{2,1}}&\longrightarrow &\boxed{H^{2,2}}&\longrightarrow &\boxed{H^{2,3}}&\longrightarrow &\cdots\\
&&&&&&&&&&\\
\boxed{H^{1,0}}&\longrightarrow &\boxed{H^{1,1}}&\longrightarrow &\boxed{H^{1,2}}&\longrightarrow &\boxed{H^{1,3}}&\longrightarrow &\boxed{H^{1,4}}&\longrightarrow &\cdots\\
&&&&&&&&&&
\end{array}$

In the diagram, vertical columns represent the breakup of $HH^n(A,M)$,
starting with $n=1$, and the horizontal arrows display the Hochschild
coboundary. The bottom row, $HH^{1,*}(A,M)$, is the Harrison
cohomology $Har^{\*}(A,M)$ which is associated to Barr's
idempotent. The idempotent $e_n(0)$ is the skew-symmetrizer
$\displaystyle \frac{1}{n!}\sum_{\sigma \in S_n} (-1)^{\sigma}\sigma$
and it follows that the diagonal components, $HH^{n,0}(A,M)$ are the
skew multi-derivations, $\bigwedge^n_A \Der (A)$, of $A$ into $M$.  If
$A=\mathcal O(V)$, the ring of regular functions on a smooth affine
variety $V$, then the celebrated Hochschild-Kostant-Rosenberg Theorem
asserts that $HH^n (A,A)=\bigwedge^n_A \Der (A)$, where $\Der(A)$. In
terms of the Hodge decomposition, the Theorem becomes
$HH^n(A,A)=HH^{n,0}(A,A)$. In particular, $Har^2(A,A)=0$ and these
algebras have no commutative deformations. In the case $V$ is not
smooth, one expects the components $HH^{r,n-r}(A,A)$ to encode
information regarding the singularities. Some interesting results by
Fronsdal in this direction can be found in \cite{Fr07}.

It is clear that the refinement of $HH^*(A,A)$ provided by the Hodge
decomposition can be useful. For example, if $HH^*(A,A)$ is
infinite dimensional, then its Euler-Poincare characteristic is not
well-defined. However, its partial Euler-Poincare characteristics
(alternating sums of $\dim H^{r,*-r}(A,A)$) may all be defined. Here
is an example which illustrates this phenomenon. Let
$A=k[\epsilon]/\epsilon^2$ be the ring of dual numbers. It is
well-known that $HH^n(A,A)$ has dimension one for all $n\geq 1$. Using
the Hodge decomposition, one can show that $HH^n(A,A)=H^{k,n-k}(A,A)$,
where $\displaystyle k=\lfloor \frac{n+1}{2}\rfloor$.  The partial
Euler-Poincare characteristics are deformation invariant and as such
they can be helpful in detecting whether a given scheme is a
deformation of another one.

N.\ Bergeron and Wolfgang showed that the components $\bigoplus_{r=1}^{k}HH^{r,n-r}(A,A)$
consist of those classes of cocycles vanishing on $(k+1)$-shuffles but not on $m$-shuffles for any $m<k+1$, see
\cite{BW95} for the precise definition and explanation. This
generalizes the fact that Harrison cohomology consists of those
cocycles vanishing on $2$-shuffles. Another fact proved in \cite{BW95}
is that $HH^{r,n-r}(A,A)$ behaves well with respect to the
filtration $\mathcal F_m=\bigoplus_{r\geq m}HH^{*,r}(A,A)$ in the
sense that $[\mathcal F_p, \mathcal F_q]\subset \mathcal F_{p+q}$.

Other instances of cohomology decompositions arising from group
actions are possible. For example, F.\ Bergeron and N.\ Bergeron found
in \cite{BB92} a type $B$ decomposition. Specifically, they showed
that there are $n$ idempotents in the descent algebra of the Weyl
group of type $B$, the group of signed permutations on $n$
letters. Moreover, if $A$ is an algebra with involution and $M$ is a
symmetric $A$-bimodule, then there is an action of $B_n$ on
$A^{\otimes n}$ with the property that the idempotents are compatible
with the Hochschild coboundary map. Thus there is a ``type $B$''
splitting of the cohomology. This raises the question of whether there
are idempotents in the descent algebras of other Coxeter systems
$(W,S)$ which decompose $HH^*(A,M)$ for algebras $A$ with a suitable $W$-action.

\section{Bialgebra Deformations}\label{bialgebra}

It was clear that, after discovery of quantum groups in the 1980's,
there should be a cohomology theory of bialgebras with the usual
features related to deformations. In \cite{GS90a} Gerstenhaber and
Schack introduced such a theory which we now describe.

The Gerstenhaber-Schack bialgebra cohomology $H^*_{GS}(B,B')$ is
defined for certain matched pairs of bialgebras $B$ and $B'$. For
simplicity, we only describe here the case $B'=B$ (any bialgebra is
matched with itself). Since $B$ is a bialgebra, any tensor power
$B^{\otimes m}$ is both a $B$-bimodule and a $B$-bicomodule and thus
the Hochschild cohomology $HH^*(B, B^{\otimes m})$ and the coalgebra
(Cartier) cohomology $H_c^*(B^{\otimes m},B)$ are well-defined. Set
$C^{p,q}(B,B)=\Hom_k (B^{\otimes p},B^{\otimes q})$.  The Hochschild
coboundary operator provides a map $\delta_h:C^{p,q}(B,B)\rightarrow
C^{p+1,q}(B,B)$ while the coalgebra coboundary yields
$\delta_c:C^{p,q}(B,B)\rightarrow C^{p,q+1}(B,B)$. These coboundaries
commute giving the Gerstenhaber-Schack complex
$$C_{GS}^{*,*}(B,B)\quad \text{with}\quad C_{GS}^n(B,B)=\mathop{\bigoplus_{p+q=n}}_{p,q>0}C^{p,q}(B,B)\quad \text{and}\quad
\delta_{GS} =\delta_h+(-1)^q\delta_c.$$
The bialgebra cohomology
$H^*_{GS}(B,B)$ is then the homology of this complex.

There are variants of this theory. For example, if one takes 
$p>0, q\geq 0$ in the definition of $C_{GS}^n(B,B)$, then the
resulting cohomology controls the deformations of $B$ to a Drinfel'd
(quasi-Hopf) algebra, see \cite{GS90b}, \cite{MS96}.  Markl has shown
in \cite{Mar07} that $H^*_{GS}(B,B)$ does carries an intrinsic graded
bracket.  In fact, Markl's construction shows the existence of a
bracket for any type of (bi)algebra over an operad or PROP.

For the rest of this section, $B$ will denote either $\mathcal O(G)$
or $U\g$, where $G$ is a reductive algebraic group and
$\g=\operatorname{Lie}(G)$. In these cases, the bialgebra cohomology
is easy to compute since $HH^n(-,\mathcal O(G))=0$ and the
$H_c^n(U\g,-)$ vanish in positive dimensions. Explicitly, if
$B=\mathcal O(G)$ or $U\g$ then
$H_{GS}^n(B,B)=\bigwedge^n\g/(\bigwedge^n \g)^{\g}$, where $(\bigwedge
\g)^{\g}$ is the space of $\g$-invariants in $\bigwedge^n\g$
\cite{GS90b}.  The Schouten bracket on $\bigwedge^*\g$ corresponds to the graded
Lie algebra structure on $H_{GS}^n(B,B)$. There are no invariants in
$\g\wedge \g$ and, up to a scalar multiple, there is a unique non-zero
invariant in $\bigwedge^3\g$. The infinitesimal bialgebra deformations
of $B$ are then elements $r\in \g\wedge \g$. The condition $[r,r]=0$
in $H^3_{GS}(B,B)$ means either that $r$ is a solution to the
classical Yang-Baxter equation (CYBE) (in the case $[r,r]=0$) or that
it is a solution to the modified CYBE (in the case that $[r,r]$ is a
non-zero invariant). Any solution to either of these Yang-Baxter
equations gives a Poisson-Lie group structure on $G$.

The quantization problem for both types of $r$-matrices is solved. For
$r$ a solution to the CYBE, the quantization is given by the UDF
associated to $r$, see Example \ref{twist} of Section \ref{UDF} and
\cite{Dr83}. For the solutions to the modified CYBE, the quantization can be
deduced from the ``dynamical twist'' found in \cite{ESS00}.  The
quantizations of \cite{Dr83} and \cite{ESS00} are universal in the
sense that they lie in $(U\g\otimes U\g)[[t]]$, and so they provide a
quantum Yang-Baxter matrix in $\operatorname{End}(V\otimes V)[[t]]$
for any representation $V$ of $\g$. Computing this $R$-matrix from the
universal quantization can require great effort. However, in
\cite{GGS93} a simple explicit ``GGS'' formula was conjectured to
quantize any modified $r$-matrix for $\g=\mathfrak{sl}(n)$ and
$V=k^n$, the vector representation. After performing computer checks
for over ten thousand cases, the GGS formula was proven correct by
Schedler in \cite{Sc00}. The proof is far from elementary as it uses
intricate combinatorial manipulations to show that the universal
solution of \cite{ESS00} coincides with the simple GGS
formula. Something is wanting for a simpler proof and real meaning of
the GGS formula. It would also be interesting to extend the result to
yield elementary quantizations of the modified $r$-matrices in the
symplectic and orthogonal cases.

The bialgebra cohomology of $\mathcal O(G)$ also guarantees that any
deformation is equivalent to one with a deformed
product $*$ which is compatible with the original comultiplication
$\Delta$. A deformation of the form $(\mathcal O(G),*,\Delta)$ is
called {\it {preferred}}. Similarly, all bialgebra deformations of
$U\g$ are preferred, although in this case it is the original
multiplication which is unchanged. The standard quantization $\mathcal
O_q(G)$ is equivalent to a preferred deformation but no such
presentation has been exhibited -- even in the simplest case of
$\mathcal O_q(\operatorname{SL}(2))$. As in the case of Lie bialgebra
quantization, the difficulty in performing explicit computations seems
to be that preferred deformations are linked with a choice of
Drinfel'd associator. See \cite{BGGS04} for a more complete discussion
of deformation quantization as it relates to quantum groups.

Returning to the Yang-Baxter equations, it should be noted that the
moduli space of solutions to the MCYBE for a simple Lie algebra has
been constructively described by Belavin and Drinfeld in
\cite{BD82}. The solutions fall into a finite disjoint union of
components, each of which is determined by an ``admissible triple''
(certain combinatorial data associated with the root system). In
contrast, an explicit classification of solutions to the CYBE is
intractable, for it would require as a special case the knowledge of
all abelian Lie subalgebras of $\g$. There is however a
non-constructive description of such $r$-matrices in terms of
``quasi-Frobenius'' Lie algebras, see \cite{BD82}, \cite{Sto91}. A Lie
algebra $\mathfrak q$ is quasi-Frobenius if there is a non-degenerate
function $\phi: \mathfrak q\wedge \mathfrak q\rightarrow k$ which is a
two-cocycle in the Chevalley-Eilenberg cohomology. The Lie algebra is
Frobenius if the two-cocycle can be taken to be a coboundary, that is, if
$\phi (a,b)=F([a,b])$ for some $F\in \mathfrak q^*$. If
$B=(B_{ij})$ is the matrix of $\phi$ with respect to some basis
$\{x_1,\dots, x_m\}$ of $\mathfrak q$, then $r=\sum
B^{-1}_{ij}x_i\wedge x_j$ is a solution to the CYBE. In \cite{GG98} it
was shown that some solutions to the CYBE arise as degenerations of
solutions to the MCYBE, and others do not. Perhaps it may be
feasible to describe all of these ``boundary'' solutions using the
Belavin-Drinfel'd triples.

\section{Diagrams of algebras}\label{diagram}

A ``diagram'' of algebras is a contravariant functor $\A$ from a small
category $\mathcal C$ to the category of associative $k$-algebras,
i.e.\ a presheaf of algebras over $\mathcal C$. So for each $i\in
\operatorname{Ob}(\C)$ there is an algebra $\A^i$ and for each
morphism $i\rightarrow j$ there is an algebra map
$\phi^{ij}:\A^j\rightarrow \A^i$. Presheaves of algebras are abundant
and surface in a variety of contexts: A single algebra is a diagram
over a one-object category with only the identity morphism. A diagram
over the category with two objects and one non-trivial morphism $u:
0\rightarrow 1$ is nothing but a homomorphism of algebras
$\phi:B\rightarrow A$. The structure sheaf $\mathcal O_{X}$ on a
quasi-projective variety $X$ is a diagram of commutative algebras over
the category $\mathcal U$ of open subsets of $X$. Here $\mathcal U$
is a category in which the morphisms correspond to inclusion maps.

In a series of papers, Gerstenhaber and Schack developed natural
cohomology and deformation theories for diagrams and proved a number
of remarkable results. A description of the theory can be found in the
survey \cite{GS88}. The Hochschild cohomology
of sheaves of algebras and abelian categories studied \cite{Hi05} and
\cite{LVdb06} are closely related to the Gerstenhaber-Schack diagram cohomology.

Perhaps the most useful and difficult result in diagram cohomology theory is
the \emph{General Cohomology Comparison Theorem} (see Theorem \ref{CCT}) which asserts, in a sense,
that the cohomology and deformation theories of an arbitrary diagram are no more general than that of a single
algebra. In order to explain more clearly what this means we give a
quick review of the basics of the theory.

An $\A$-bimodule $\M$ is a contravariant functor from $\C$ to the
category of abelian groups assigning to every $i\in
\operatorname{Ob}(\mathcal C)$ an $\A^i$-bimodule $\M^i$ and to every
morphism $u: i\rightarrow j$ in $\mathcal C$ a map
$T^{ij}:\M^j\rightarrow \M^i$ which is required to be an
$\A^j$-bimodule map. Here, $\M^i$ becomes an $\A^j$-module by virtue
of the algebra homomorphism $\phi^{ij}$.

Just as in the case of a single algebra, there are various descriptions
of the diagram cohomology $H_d^*(\A,\M)$. Once the requisite
categorical machinery is laid out, one may define
$H_d^*(\A,\M)=\operatorname{Ext}_{\A-\A}(\A,\M)$. There is also a
cochain description which is quite useful and we present this
here. There is a cochain complex $(C_d^*(\A,\M),\delta_d)$ whose
homology coincides with $\operatorname{Ext}_{\A-\A}(\A,\M)$. The
description of $C_d^*(\A,\M)$ has both algebraic and simplicial
aspects. The nerve $\Sigma$ of $\mathcal C$ is the simplicial complex
whose $0$-simplices are the objects of $\mathcal C$ and the
$p$-simplices are the composable maps $\sigma =(i_0\rightarrow
i_i\rightarrow \cdots \rightarrow i_p)$. For simplicity we write
$\sigma =(i_0,\ldots, i_p)$.  The boundary of $\sigma$ is $\partial
\sigma=\sum(-1)^j\sigma_j$, where $\sigma_j$ is the $j$-th face of
$\sigma$ obtained by omitting $i_j$.

For a diagram $\A$ and $\A$-bimodule $\M$, the $n$-cochains are
$C_d^n(\A,\M)=\bigoplus_{p+q=n}C^{p,q}_d(\A,\M),$ where

$$C_d^{p,q}(\A,\M)=\prod_{\substack {p-\text{simplices}\\(i_0,\ldots,i_p)}}C^q(\A(i_p),\M(i_0)).$$

Fix $\Gamma \in C_d^{p,q}(\A,\M)$. The diagram coboundary will have
two components: $\delta_{alg}\Gamma\in C_d^{p,q+1}(\A,\M)$ and
$\delta_{simp}\in C_d^{p+1,q}(\A,\M)$.  The algebraic component is
defined by $(\delta_{alg}\Gamma)^{\sigma}=\delta_h(\Gamma^{\sigma})$
where $\delta_h:C^q(\A^{i_p},\M^{i_0})\rightarrow
C^{q+1}(\A^{i_p},\M^{i_0})$ is the ordinary Hochschild coboundary
operator. The simplicial component is defined as follows. Let
$\sigma=(i_0,\ldots,i_{p+1})$ be a $p+1$-simplex. For faces $\sigma_j$
with $1\leq j\leq p$, we have $\Gamma^{\sigma_j}\in
C^q(\A^{i_{p+1}},\M^{i_0})$, while $\Gamma^{\sigma_0}\in
C^q(\A^{i_{p+1}},\M^{i_1})$ and $\Gamma^{\sigma_{p+1}}\in
C^q(\A^{i_{p}},\M^{i_0})$. The extreme cases $\Gamma^{\sigma_0}$ and 
$\Gamma^{\sigma_{p+1}}$ lie in different cochain groups than the others, but
there are adjustments however which correct
this. For $\sigma_0$ note that the composite
$T^{i_0i_1}\Gamma^{\sigma_0}\in C^{q}(\A^{i_{p+1}},\M^{i_0})$. For
$\sigma_{p+1}$ define $\Gamma^{\sigma_{p+1}}\phi^{i_{p+1}i_p}\in
C^q(\A^{i_{p+1}},\M^{i_0})$ by
$$\Gamma^{\sigma_{p+1}}\phi^{i_{p+1}i_p}(a_1,\ldots,a_q)=\Gamma^{\sigma_{p+1}}(\phi^{i_{p+1}i_p}a_1,\ldots,
\phi^{i_{p+1}i_p}a_q).$$ Now set
$$(\delta_{simp}\Gamma)^{\sigma}=\text{``}\Gamma^{\partial
  \sigma}\text{''}=T^{i_0i_1}\Gamma^{\sigma_0}-\Gamma^{\sigma_1}+\Gamma^{\sigma_2}-\cdots
+(-1)^p\Gamma^{\sigma_p}\phi^{i_{p}i_p+1}.$$ The full diagram
coboundary is then
$$\delta_d = \delta_{alg}+(-1)^p\delta_{simp}$$
and the diagram cohomology $H_d^*(\A,\M)$ is defined to be the homology of the complex
$$C_d^*(\A,\M)=\bigoplus _{p+q=n}C_d^{p,q}(\A,\M)\quad \mbox{with}
\quad \delta_d =\delta_{alg}+(-1)^p\delta_{simp}.$$ Note that the
cohomology of the bottom row $H_d^{*,0}(\A,\M)$ coincides with the
simplicial cohomology of $\Sigma(\C)$ with local coefficients $\M$.

A deformation of $\A$ is a diagram of $k[[t]]$-algebras whose
reduction modulo $t$ is $\mathbb{A}$. The diagram cohomology
$H^*_d(\A,\A)$ is too large to govern deformations of $\A$ since the
simplicial cohomology of $\Sigma(\mathcal C)$ may not be trivial. There are
remedies such as using ``asimplicial'' cochains or adjoining a
terminator to $\mathcal C$, see \cite{GS88}. Naturally, we would like
a graded Lie structure on $H^*_d(\A,\A)$ which controls
obstructions. It turns out that, unlike the case of a single algebra,
the natural bracket on $C^*_d(\A,\A)$ gives the structure of only a
homotopy graded Lie algebra. Proving that this bracket descends to a
graded Lie structure at the cohomology level would be at best a nasty
computation using the cochain description. However, the following very
difficult and useful result of \cite{GS88} settles this question.

\begin{theorem}[General Cohomology Comparison Theorem]\label{CCT} Associated to each diagram $\A$ is a
single $k$-algebra $\A!!$ such that the cohomology and deformation
theories of $\A$ are naturally isomorphic to those of $\A!!$. In
particular, $H^*_d(\A,\A)$ is a Gerstenhaber algebra.
\end{theorem}

The \textit{diagram algebra} $\A!!$ is rather complicated and we will
not describe it here, although we will see a special case in Section
\ref{variation}. The proof of Theorem \ref{CCT} relies on the
\textit{Special Cohomology Comparison Theorem} which is the case when $\A$ is a poset. To derive
the general case, Gerstenhaber and Schack perform a barycentric subdivision of $\A$. It turns out that
the second subdivision of an arbitrary diagram is a poset and subdivision preserves the cohomology. Van den Bergh and Lowen
have proved Special Cohomology Comparison Theory for prestacks in \cite{LVdb09}.

Another important result in diagram cohomology theory is the following theorem
which completely reconciles the Kodaira-Spencer manifold deformation
theory with the Gerstenhaber-Schack diagram deformation theory.

\begin{theorem} Let $X$ be a smooth compact algebraic variety with tangent bundle $T$.
Suppose $\mathcal U$ be an affine open cover of $X$ and let $\A$ be
the restriction of $\mathcal O_X$ to $\mathcal U$. Then there is a
Gerstenhaber algebra isomorphism $H^*_d(\A,\A)\simeq H^*(X,\bigwedge^*
T)$.
\end{theorem}

Using the theorem, one sees that
$$H^2_d(\A,\A)\simeq H^2(X,\mathcal O_X)\bigoplus H^1(X,T)\bigoplus H^0(X,\wedge ^2T).$$

The middle term consists of the infinitesimal deformations of $X$ in
the Kodaira-Spencer theory. The last term is the space of
infinitesimal deformations of $X$ to ``non-commutative'' spaces; those
with vanishing primary obstruction are precisely the Poisson
structures on $X$ and, by Theorem \ref{Kontsevich}, these are
quantizable. The meaning of the first term of $H^2_d(\A,\A)$ is not
well-understood.

Besides applications to geometric situations, diagrams naturally arise
in other contexts. For example, given an algebra $A$ and an $A$-module $M$,
one can deform the action of $A$ on $M$ in the evident way, and it is
relatively easy in this case to deduce the appropriate deformation
cohomology. More generally, one can simultaneously deform $A$ and its
action on $M$ in a compatible way. These situations are special cases
of diagram deformations. Indeed, the original $A$-module structure on
$M$ is simply an algebra homomorphism $\phi:A\rightarrow \End(M)$, and
hence is a diagram. Deformations of this diagram yield the various
possibilities of deforming $A$, the action of $A$ on $M$, or both. The
general theory automatically yields appropriate cohomology and
obstruction theories. In Section \ref{variation} the diagram cohomology
theory will be used to cohomologically explain how certain rigid
algebras can appear in naturally parametrized families.

\section{Deforming relations}\label{relations}

Suppose an algebra is given as $A=TX/J$ where $X$ is the $k$-module
spanned by finitely many generators $x_i$, $TX$ is the tensor algebra,
and $J$ is the ideal of relations. If $J_t$ is an ideal of $TX[[t]]$
which reduces to $J$ modulo $t$, then a natural question is whether
$A_t=TX[[t]]/J_t$ is a deformation of $A$ or not. Associativity of
$A_t$ is automatic but to be a deformation it must be flat as a
$k[[t]]$-module. There is no efficient way in general to determine
if the relations in $J_t$ insure flatness.  An elementary
case where flatness fails is the following: Let $A=k[x,y,z]$ and let
$J_t$ be generated by $yx-(1+t)txy$, $zx-xz-ty^2$ and $yz-zy$. When
$t=0$ all variables commute and the polynomial algebra $k[x,y,z]$ is
obtained. For $t\neq 0$, the deformed relations allow for a PBW-type ordering in
which every monomial of $A_t$ can be reduced to one of the form
$x^iy^kz^k$. However, the element $t(1+t)y^3$ lies in $J_t$ and so
$A_t$ has $t$-torsion and thus is not flat.

Flatness is relatively easy to check for certain deformations of
Koszul algebras, which comprise an important class of quadratic
algebras. An algebra $A$ is quadratic if $A=TX/J$, with $J$ is
generated by relations $R\subset X\otimes X$. Since the relations are
homogeneous, such algebras are $\N$-graded, $A=\bigoplus_{i\geq
  0}A[i]$ and $\dim A[i]<\infty $ for each $i$. In particular,
$A[0]=k$. A quadratic algebra $A$ is \emph{Koszul} if its dual $A^!$
is isomorphic to the Yoneda algebra $\operatorname{Ext}^*_A(k,k)$.
Variations of the following fundamental theorem have appeared in
several places in the literature, most notably in the works of
Drinfel'd \cite{Dr86} and Braverman-Gaitsgory \cite{BG96}.

\begin{theorem}[Koszul Deformation Criterion]\label{koszul}
Suppose that $A=TX/J$ is Koszul and $A_t=TX[[t]]/J_t$, where $J_t$ is
generated by relations $R_t\subset (X\otimes X)[[t]]$ which reduce to
$J$ modulo $t$. Then $A_t$ is a deformation of $A$ if and only if
$A_t[3]$ is a flat $k[[t]]$-module.
\end{theorem}
The point of the theorem is that in the Koszul case, flatness in
dimensions greater than 3 is a consequence of flatness in degree 3. Flatness
in the cases of degrees 1 and 2 is automatic.

One of the most interesting and explicit uses of the Koszul
deformation criterion has been carried out by Etingof and Ginzburg in
the theory of \emph{symplectic reflection algebras}, which are
deformations of crossed product algebras $SV\ltimes G$, see
\cite{EG02}.  One can try to deform $SV\ltimes G$ by imposing
additional relations of the form $xy-yx = \kappa(x,y)$ where $x,y\in V$ and
$\kappa(x,y) = -\kappa(y,x)\in \mathbb CG$. For an arbitrary
skew-symmetric function $\kappa$, the underlying vector space of the
resulting algebra, $A_{\kappa}$, will be smaller than that of $A_0=
SV\ltimes G$ -- that is, the deformation will not be flat.

In the case where $V$ is a symplectic vector space and $G\in Sp(V)$,
Etingof and Ginzburg have an explicit and remarkable classification of
which skew forms $\kappa$ lead to deformations. To describe these, we
first need some notation.  Suppose $V$ is a complex vector space
equipped with a skew bilinear form $\omega:V\times V\rightarrow
\mathbb C$, and let $G$ be a finite subgroup of $Sp(V)$. An element
$s\in G$ is a \emph{symplectic reflection} if the rank of $1-s$ is
2. The set of all symplectic reflections is denoted $S$. For each
$s\in S$, let $\omega_s$ denote the form on $V$ with radical
$\operatorname{Ker}(1-s)$ and which coincides with $\omega$ on
$\operatorname{Im}(1-s)$. The triple $(V,\omega, G)$ is
\emph{indecomposable} if $V$ can not be split into a non-trivial
direct sum of $G$-invariant symplectic subspaces.

\begin{theorem}[Etingof-Ginzburg] Suppose $(V,\omega, G)$ is an
  indecomposable triple, and let $\kappa: V\times V\rightarrow \mathbb
  C G$ be a skew form. Then $A_{\kappa}$ is a flat deformation of
  $SV\ltimes G$ if and only if there exists a $G$-invariant function
  $c:S\rightarrow \mathbb C$, $s\mapsto c_s$ and a constant $t$, such
  that
$$\kappa(x,y)=t\omega(x,y) +\sum_{s\in S}c_s\omega_s(x,y)s.$$
\end{theorem}

As stated earlier, the applications of symplectic reflection algebras
are many. Here is one particularly interesting one. The center of
$SV\ltimes G$ is the algebra $(SV)^G$ of $G$-invariant polynomial
functions, which can be viewed as the functions on the orbit space
$V/G$. If $e=\frac{1}{|G|}\sum_{g\in G}g$ is the symmetrizing
idempotent in $\mathbb CG$, then the \emph{spherical} subalgebra of
$A_{\kappa}$ is defined to be $eA_{\kappa}e$. It is known that $eA_0e\simeq (SV)^G$,
and so $eA_{\kappa}e$ provides a non-commutative deformation of
$(SV)^G$. However, if $t=0$ then the algebra $eA_{\kappa}e$ is
commutative. Thus the symplectic reflection algebras can provide
geometric deformations of $V/G$.

Returning to Theorem \ref{koszul}, there are algebras
where, unlike the symplectic reflection algebras, there is no evident
ordered or PBW-type basis of $A_t$. For example, Sklyanin (or
elliptic) deformations of polynomial algebras have this property. The
simplest case is the algebra with generators $\{x,y,z\}$ and relations
$$ax^2+byz+czy=0,\qquad ay^2+bzx+cxz=0,\qquad az^2+bxy+cyx=0.$$ The
triple $(a,b,c)=(0,1,-1)$ gives the polynomial algebra $k[x,y,z]$, but
for generic $(a,b,c)$ the relations are such that there is no PBW-type
basis.  One way to prove flatness is to associate certain geometric
data (an elliptic curve $\mathcal E$ and point $\eta\in \mathcal E$)
to the algebra in question. The geometric information allows one to
construct a factor ring of the Sklyanin algebra which can be exploited
to establish flatness. A survey of elliptic deformations of polynomial
algebras can be found in \cite{Od02}.

\section{Variation of algebras}\label{variation}

As mentioned in section \ref{rigid}, an algebra with $H^2(A,A)=0$ may
depend essentially on parameters and so the classic deformation theory
of $A$ does not detect this dependence. If we instead pass to an
appropriate diagram of algebras, it is possible in many cases to
detect the dependence of $A$ on parameters from the diagram cohomology
and construct the new algebra with the concept of algebra variation.
In this section we give a brief account of \cite{GG08b}.

Suppose that we have $k$-algebras $A, B, B'$ and monomorphisms $\phi:B
\to A$ and $\phi':B' \to A$ such that $A$ is generated by the images
$\phi(B)$ and $\phi'(B')$. If $V$ is the direct sum of the underlying
$k$-modules of $B$ and $B'$ then $A=TV/J$, where $J$ is the ideal of
$TV$ generated by relations which we write in the form $R(\phi(b),\phi(b'))$ for
$b\in B$ and $b'\in B'$. In this case we have a diagram $\mathbb A$
over the poset $\mathcal C=\{0,1,1'\}$:

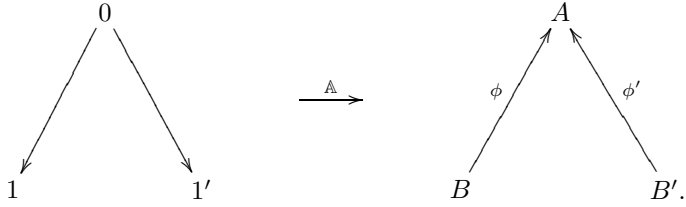
$$\xymatrix{& 0\ar[ldd] \ar[rdd]&&&&&A&\\ &&&\ar[r]^{\mathbb
    A}&&&&\\ 1&& 1'&&&B\ar[ruu]^{\phi}&& \ar[luu]_{\phi'}B'.}\label{poset}$$

Now consider a deformation $\mathbb A_t$ of $\mathbb A$ in which the
algebras $A$, $B$, and $B'$ remain fixed but the homomorphism $\phi$ is
deformed as $\phi_t=\phi+t\phi_1+t^2\phi_2+\cdots$ and similarly
assume $\phi'$ is deformed to $\phi'_t$. We can use the \textit{same}
relations determining $A$ with deformed inputs  to construct a new algebra $A_t$.

\begin{definition} Suppose $A, B, B', V, R, \mathbb A$ and $ \mathbb A_t$ are as above. Let
$J_t$ be the ideal of $TV[[t]]$ generated by all elements of the form
  $R(\phi_t(b),\phi_t'(b'))$ for $b\in B$ and $b'\in B'$. The algebra
  $A_t=TV[[t]]/J_t$ is called a \textrm{variation} of $A$.
\end{definition}

A variation $A_{t}$ is certainly associative but there is no guarantee
that it is flat, and as noted earlier, there is in general no easy way
to determine when such algebras are flat. The concept of variation
can clearly be generalized by letting $A$ be generated by more than two
subalgebras.

It is important to note that not all algebras of the form $TV[[t]]/J_t$
where $J_t$ is an ideal of $TV[[t]]$ with $J_0=J$ are variations of
$A$. As an example, take $A$ to be commutative. Then we have in $J$
all relations of the form $\phi(b)\phi'(b')-\phi'(b')\phi(b)$.  The
ideal $J_t$ defining the variation $A_t$ will therefore have all
relations of the form $\phi_t(b)\phi_t'(b')-\phi_t'(b')\phi_t(b)$ and
so $A_t$ remains commutative.

Let us return now to the deformation of the diagram $\A$ in the above figure obtained by
replacing $\phi$ with $\phi_t$ and $\phi'$ with $\phi_t'$. Its infinitesimal
lies in $ H^2_d(\mathbb A, \mathbb A)$ and is the class of a cocycle of the form
$\Gamma =(\Gamma_A, \Gamma_B, \Gamma_{B'}, \Gamma_{BA}, \Gamma_{B'A})$ with
\begin{multline} \Gamma_A\in HH^2(A,A), \quad \Gamma_B\in HH^2(B,B), \quad
\Gamma_{B'}\in HH^2(B',B'), \\ \Gamma_{BA}\in HH^1(B,A)
\quad\mbox{and}\quad \Gamma_{B'A}\in HH^1(B',A).\end{multline}
The first three
components of $\Gamma$ have algebraic dimension $2$ and simplicial
dimension $0$ while the last two have algebraic and simplicial
dimension $1$ as these correspond to the $1$-simplices of the
underlying category. The deformation $\A_t$ may be viewed as an \textit{integral}
of this cohomology class. We also assign this class to the variation $A_t$.

Even if the algebras $A$, $B$, and $B'$ are absolutely rigid,
$H^2_d(\A,\A)$ may not vanish in general as
$HH^1(B,A)=\operatorname{Der}(B,A)\neq 0$ and similarly for
$HH^1(B',A)$. In this case, $\Gamma$ obviously can be taken to
be of the form $(0,0,0, \Gamma_{BA},\Gamma_{B'A}).$ However, if the
characteristic of $k$ is zero, then we may further assume that $\Gamma
=(0,0,0,0,\Gamma_{B'A})$.

\begin{remark}The diagram algebra $\A !!$ associated to $\A$ (see
  Theorem \ref{CCT}) can be viewed as the algebra of $3\times 3$
  matrices of the form
$$\left[
\begin {array}{ccc}
a_1&a_2&a_3\\
\noalign{\medskip}
0&b&0\\
\noalign{\medskip}
0&0&b'
\end {array}
\right], \quad \mbox{with}\quad a_i\in A,\,\, b\in B,\,\,b'\in B',$$
where the multiplication in $\A !!$ uses the
convention that $ba=\phi(b)a$ and similarly for $b'a$.  Even in this
simple case it is difficult to see how to canonically relate the
cohomology and deformations of $\A$ with those of $\A !!$.
\end{remark}

We end with a reconsideration of the first quantized Weyl algebra
$A_q=\C[x,y]$ with relation $qxy-yx=1$. Having already remarked that
$A_q$ is not isomorphic to $A_1$ for $q$ near 1, yet
$HH^2(A_1,A_1)=0$. Consider now whether $A_1$ can be varied to
$A_q$. Using our earlier notation, suppose $A=A_1$, $B=\C[x]$, and
$B'=\C[y]$ and let $\phi:B\rightarrow A$ and $\phi:B'\rightarrow A$ be
the inclusion maps. All of these algebras are absolutely rigid. Thus,
based on the comments above, it suffices to vary the inclusion morphism of
$\mathbb{C}[y]$ into $A_1$. The question becomes whether there exists
an element $y'\in A_1[[t]]$ of the form $y+t\eta_1 +t^2\eta_2 +
\cdots, \eta_i \in A_1$ such that the relation $[x, y']= xy'-y'x
= 1$ is equivalent to having $[x,y] = 1-t xy$, for this would give
$(1+t)xy - yx = 1$, i.e., $A_q[[t]]$ with $q=1+t$. There are indeed
elements $y'$ of the desired form.  In \cite{GG08b}, it is shown that
one may take
$$y' = y + a_1(t)xy^2 + a_2(t)x^2y^3 + \dots, \quad \mathrm{\
where\ \ } a_r(t) = \frac{t^{r+1}}{(1+t)^{r+1}-1}.$$

Thus, $A_q$ is a variation of $A_1$. It is easy to use the formula for
$y'$ to show that the corresponding diagram infinitesimal is $\Gamma =
(0,0,0,0,\delta)$, where $\delta\in \operatorname{Der}(\C[y],A_1)$ is
the derivation with $\delta(y)= xy^2$. This is a non-trivial cohomology class in
$H^2_d(\A,\A)$ and so the diagram cohomology has detected the variation from $A_1$ to $A_q$.

It is instructive to note that
in the power series representation of the
$a_r(t)$ there could be no value of $t$ for which all the
series converge, for each $a_r(t)$ is a rational function with
a pole wherever $t$ has the form $\omega - 1$ where $\omega$
is an $(r+1)$st root of unity and every neighborhood of 0 in
$\mathbb{C}$ contains infinitely many of these. Those $A_q$ with
$q$ a root of unity are in some sense `unreachable' from $A_1$.
Nevertheless, $y'$ can actually be evaluated for any complex number
$t$ with $1+t$ not a root of unity.


\begin{thebibliography}{MMXS06}

\bibitem[Ba68]{Ba68} Barr M.: Harrison homology, Hochschild homology,
  and triples, J. Alg. \textbf{8}, (1968), 314-323.

\bibitem[BFFLS78]{BFFLS78} Bayen F., Flato M., Fronsdal C., Lichnerowicz
  A., Sternheimer D.: {Deformation theory and quantization I:
    Deformation of symplectic structures; II: Physical applications},
  Ann. Phys. (NY) \textbf{111} (1978), 61--110; 111--151.



\bibitem[BD82]{BD82} Belavin A. and Drinfel{'}d V., Solutions of the
  classical Yang-Baxter equation for simple Lie algebras. (Russian),
  \textit{Funktsional. Anal. i Prilozhen.} \textbf{16} (1982) no. 3,
  1--29.

\bibitem[BB92]{BB92} F. Bergeron and N. Bergeron.: Orthogonal
  idempotents in the descent algebra of type $B_n$ and applications,
  J. Pure Appl. Alg., \textbf{79}(2) (1992), 109-129.

\bibitem[BW95]{BW95} Bergeron N. and Wolfgang H.L.: The decomposition
  of Hochschild cohomology and Gerstenhaber operations, J. Pure
  Appl. Algebra, \textbf{104}(3) (1995), 243--265.

\bibitem[BTY07]{BTY07} Bieliavsky, P., Tang, X., and Yao, Y.:
  Rankin-Cohen brackets and formal quantization.
  Adv. Math. \textbf{212} (2007), no. 1, 293--314.



\bibitem[BGGS04]{BGGS04} Bonneau P., Gerstenhaber M., Giaquinto A.,
  Sternheimer D.: Quantum groups and deformation quantization:
  explicit approaches and implicit aspects.  J. Math. Phys.
  \textbf{45}(10) (2004), 3703--3741.


\bibitem[BG96]{BG96} Braverman A. and Gaitsgory D.:
  Poincare-Birkhoff-Witt theorem for quadratic algebras of Koszul
  type.  J. Algebra \textbf{181} (1996), no. 2, 315--328.



\bibitem[CGW04]{CGW04} Caldararu A., Giaquinto A., and Witherspoon S.:
  Algebraic deformations arising from orbifolds with discrete torsion, 
  J. Pure Appl. Alg., \textbf{187}(1) (2004), 51-70.

\bibitem[CF00]{CF00} Cattaneo A.S. and Felder G.: A path integral
  approach to the Kontsevich quantization formula, Comm. Math. Phys.,
  \textbf{212} (3) (2000) 591--611.

\bibitem[CKTB05]{CKTB05} Cattaneo A., Keller B., Torossian C., and
  Bruguieres A.: D\'{e}formation, quantification, th\'{e}orie de Lie. (French)
  [Deformation, quantization, Lie theory] Panoramas et Synth\`{e}ses
  [Panoramas and Syntheses], 20. Soci\'{e}t\'{e} Math\'{e}matique de France,
  Paris, 2005.

\bibitem[CM04]{CM04} Connes A. and Moscovici H.: Rankin-Cohen brackets
  and the Hopf algebra of transverse geometry,
  Mosc. Math. J. \textbf{4} (2004), no. 1, 111--130.


\bibitem[Do60]{Do60} Douady A.: Obstruction primaire \'a la
  d\'eformation, S\'eminaire Henri Cartan, \textbf{ 4}, 1960-61,
  expos\'e.

\bibitem[DMZ07]{DMZ07} Doubek, M., Markl, M., and  Zima, P.:
Deformation theory (lecture notes),
Arch. Math. (Brno) \textbf{43} (2007), no. 5, 333--371.


\bibitem[Dr83]{Dr83} Drinfel'd V.\ G.: On constant, quasiclassical
  solutions of the classical Yang-Baxter equation, Soviet Math. Dokl.,
  \textbf{28}(3), 667-671.

\bibitem[Dr86]{Dr86} Drinfel'd V.\ G.: On quadratic commutation
  relations in the quasiclassical case [translation of Mathematical
    physics, functional analysis (Russian), 25--34, 143, ``Naukova
    Dumka'', Kiev, 1986; MR0906075 (89c:58048)]. Selected
  translations.  Selecta Math. Soviet.\ \textbf{11} (1992), no. 4, 317--326.

\bibitem[Dr92]{Dr92} Drinfe'd V.G.: On some unsolved problems in
  quantum group theory, Lect. Notes Math., \textbf{1510} (1992), 1-8.


\bibitem[EG02]{EG02} Etingof P. and Ginzburg V.: Symplectic reflection
  algebras, Calogero-Moser space, and deformed Harish-Chandra
  homomorphism, Invent. Math., \textbf{147} (2002), no. 2, 243--348.

\bibitem[EK96]{EK96} Etingof P. and Kazhdan D.: Quantization of Lie
  bialgebras. I., Selecta Math. (N.S.) \textbf{2} (1996),
  no. 1, 1--41.

\bibitem[ESS00]{ESS00} Etingof P., Schedler T. and Schiffman O.:
  Explicit quantization of dynamical $r$-matrices for finite
  dimensional semisimple Lie algebras,
  J. Amer. Math. Soc. \textbf{13}(3) (2000), 595--609.


\bibitem[FN57]{FN57} Fr\"olicher A. and Nijenhuis A.:
A theorem on stability of complex structures.
Proc. Nat. Acad. Sci. \textbf{43} (1957), 239--241.


\bibitem[Fr07]{Fr07} Fronsdal C.: Quantization on curves,
  Lett. Math. Phys. \textbf{79} (2007), no. 2, 109--129. (appendix by
  Kontsevich M.)

\bibitem[Ga90]{Ga90} Garsia A.: Combinatorics of the free Lie algebra
  and the symmetric group, \textit{Analysis, et Cetera}, 309--382,
  Academic Press, Boston, MA, 1990.


\bibitem[Ge63]{Ge63} Gerstenhaber M.: The cohomology structure of an
  associative ring, Ann. Math. \textbf{78} (1963), 267--288.

\bibitem[Ge64]{Ge64} Gerstenhaber M.: On the deformation of rings and
  algebras, Ann. Math. (2) \textbf{79} (1964), 59--103.


\bibitem[Ge68]{Ge68} Gerstenhaber M.: On the deformation of rings and
  algebras III, Ann. Math. (2) \textbf{88}(1), 1-34.


\bibitem[GG98]{GG98} Gerstenhaber M. and Giaquinto A,: Boundary
  solutions of the quantum Yang-Baxter equation and solutions in three
  dimensions, \textit{Lett. Math. Phys.} \textbf{44} (1998), no. 2,
  131--141.

\bibitem[GG08a]{GG08a} Gerstenhaber M. and Giaquinto A.: Graphs,
  Frobenius functionals, and the classical Yang-Baxter equation,
  \texttt{arXiv:0808.2423v1}.


\bibitem[GG08b]{GG08b} Gerstenhaber M. and Giaquinto A.: Variation of
  Algebras, preprint, 2008.

\bibitem[GGS93]{GGS93} Gerstenhaber M., Giaquinto A. and Schack S.D.:
  Construction of quantum groups from Belavin-Drinfel{'}d
  infinitesimals, \textsl{Quantum deformations of algebras and their
    representations}, Israel Math. Conf. Proc. \textbf{7}, 45-64,
  Bar-Ilan Univ., Ramat Gan, 1993.


\bibitem[GS86]{GS86} Gerstenhaber, M. and Schack, S. D.: Relative
  Hochschild cohomology, rigid algebras and the Bockstein, J. Pure
  Appl. Alg. \textbf{43} (1986), 53-74.

\bibitem[GS87]{GS87} Gerstenhaber, M. and Schack, S. D.: A Hodge-type
  decomposition for commutative algebra cohomology, J. Pure Appl. Alg.
  \textbf{48} (1987), 229-247.

\bibitem[GS88]{GS88} Gerstenhaber M., Schack S.D.: Algebraic
  cohomology and deformation theory.  \textit{Deformation theory of
    algebras and structures and applications} (Il Ciocco, 1986),
  11--264, NATO Adv. Sci. Inst. Ser. C Math. Phys. Sci.  \textbf{247},
  Kluwer Acad. Publ. Dordrecht 1988.

\bibitem[GS90a]{GS90a} Gerstenhaber M., Schack S.D.: Bialgebra
  cohomology, deformations, and quantum groups, Proc. Nat. Acad. Sci,
  \textbf{87}(1) (1990), 478-481.


\bibitem[GS90b]{GS90b} Gerstenhaber M., Schack S.D.: Algebras,
  bialgebras, quantum groups, and algebraic deformations,
  \textit{Deformation theory and quantum groups with applications to
    mathematical physics (Amherst, MA, 1990)}, 51--92, Contemp. Math.,
  \textbf{134}, Amer. Math. Soc., Providence, RI, 1992.


\bibitem[GZ98]{GZ98} Giaquinto A. and Zhang J.: Bialgebra actions,
  twists, and universal deformation formulas, Jour. Pure
  Appl. Alg. \textbf{128}(02) (1998), 133-151.


\bibitem[Hi05]{Hi05} Hinich V.: Deformations of sheaves of algebras,
  Adv. Math. \textbf{195}(1) (2005), 102--164
  (\texttt{math.AG/0310116}).


\bibitem[KS58]{KS58} Kodaira K. and Spencer D.C.: On deformations of
  complex analytic structures, {Ann. Math.} (2) {\bf 67} (1958),
  328--466.



\bibitem[Ko97]{Ko97} Kontsevich M.: Deformation quantization of
  Poisson manifolds, {Lett. Math. Phys.} \textbf{66}(3) (2003),
  157--216 (\texttt{q-alg/9709040}).


\bibitem[Ko01]{Ko01} Kontsevich M.: Deformation quantization of
  algebraic varieties, EuroConf\'erence Mosh\'e Flato 2000, Part III
  (Dijon), Lett. Math. Phys.  \textbf{56}(3) (2001), 271--294.




\bibitem[KLO01]{KLO01} Kulish P. P., Lyakhovsky V.D. and del Olmo
  M. A.: Chains of twists for classical Lie algebras, J. Phys. A:
  Math. Gen. \textbf{32} (1999), 8671--8684.

\bibitem[LVdb06]{LVdb06} Lowen W. and Van den Bergh M.: Deformation
  theory of abelian categories, Trans. Amer. Math. Soc. \textbf{358}
  (2006), 5441-5483.

\bibitem[LVdb09]{LVdb09} Lowen W. and Van den Bergh M.: A Hochschild
  cohomology Comparison Theorem for prestacks,
  Trans. Amer. Math. Soc., to appear.



\bibitem[LS02]{LS02} Lyakhovsky V. D. and Samsonov M. E.:
Elementary parabolic twist, J. Algebra Appl. \textbf{1} (2002), no. 4, 413--424.

\bibitem[Lo88]{Lo88} Loday J-L.: Partition eul\'{e}rienne et
  op\'{e}rations en homologie cyclique, C. R. Acad. Sci. Paris S\'{e}r. I
  Math. \textbf{307}(7) (1988), 283--286.

\bibitem[Lo89]{Lo89} Loday J-L.: Op\'{e}rations sur l'homologie
  cyclique des alg\`{e}bres commutatives, Invent. Math. \textbf{96}(1)
  (1989), 205-230.



\bibitem[Mat97]{Mat97} Mathieu O.: Homologies associated with Poisson
  structures, \textit{Deformation theory and symplectic geometry
    (Ascona, 1996)}, Math. Phys. Stud., \textbf{20}, Kluwer
  Acad. Publ., Dordrecht, (1997), 177-199.

\bibitem[Mar07]{Mar07} Markl M.: Intrinsic brackets and the
  $L_\infty$-deformation theory of bialgebras, arXiv:math/0411456v6


\bibitem[MS96]{MS96} Markl M. and Shnider S.: Cohomology of Drinfeld
  algebras: a homological algebra approach,
  Internat. Math. Res. Notices (9) (1996), 431--445.

\bibitem[Od02]{Od02} Odesskii, A. V., Elliptic algebras. (Russian)  
Uspekhi Mat. Nauk  \textbf{57}  (2002),  no. 6(348), 87--122;  
translation in  Russian Math. Surveys  \textbf{57}  (2002),  no. 6, 1127--1162.


\bibitem[Pi06]{Pi06} Pinczon G.: On two theorems about symplectic reflection algebras,
\texttt{math.QA/0612690}.

\bibitem[Ri67]{Ri67} Richardson R.W.: On the rigidity of semi-direct
  products of Lie algebras, Pacific J. Math., \textbf{22}(2) (1967),
  339-344.

\bibitem[Ri57]{Ri57} Riemann B.: Theorie der Abel'schen Functionen,
  Journal f\"ur die reine und angewandte Mathematik, \textbf{54} (1857),
  101-155)



\bibitem[Sc00]{Sc00} Schedler T., Proof of the GGS
  conjecture. \textit{Math. Res.  Lett.} \textbf{7} (2000), no. 5-6,
  801--826.

\bibitem[Ste]{Ste} Sternheimer D.: The deformation philosophy, quantization
and noncommutative space-time structures, this volume.


\bibitem[Sta63]{Sta63} Stasheff J.,:
Homotopy associativity of $H$-spaces. I, II.
Trans. Amer. Math. Soc. \textbf{108} (1963), 275-292;



\bibitem[Sta93]{Sta93} Stasheff J.: The intrinsic bracket on the
  deformation complex of an associative algebra, J. Pure Appl. Algebra
  \textbf{89} (1993), 231–235.

\bibitem[Sto91]{Sto91} Stolin A., On rational solutions of Yang-Baxter
  equation for $\mathfrak{sl}(n)$, \textit{Math. Scand.} \textbf{69}
  (1991), 57--80.



\bibitem[VdB07]{VdB07} Van den Bergh M.: On global deformation
  quantization in the algebraic case, J. Algebra \textbf{315}(1)
  (2007), 326--395 (\texttt{math.AG/0603200}).


\bibitem[Ye05]{Ye05} Yekutieli A.: Deformation quantization in
  algebraic geometry,  Adv. Math.  \textbf{198}(1) (2005), 383--432.


\end{thebibliography}
\end{document}